\documentclass[12pt,leqno,twoside]{amsart}
\usepackage{amssymb}
\usepackage{latexsym}
\usepackage{verbatim}
\usepackage{epsfig}

\newtheorem{theorem}{Theorem}[section]
\newtheorem{corollary}[theorem]{Corollary}
\newtheorem{lemma}[theorem]{Lemma}
\newtheorem{proposition}[theorem]{Proposition}
\theoremstyle{definition}
\newtheorem{definition}[theorem]{Definition}
\theoremstyle{remark}
\newtheorem{remark}[theorem]{\sc Remark}
\theoremstyle{remark}
\newtheorem{example}[theorem]{\sc Example}
\theoremstyle{remark}

\theoremstyle{remark}

\theoremstyle{remark}

\numberwithin{equation}{section}

\renewcommand{\Box}{\square}    

\newcommand{\cal}{\mathcal}

\renewcommand{\int}{{\rm{int}}}

\newcommand{\pol}{{\rm{pol}}}

\newcommand{\Sing}{{\rm{Sing\hspace{2pt}}}}

\newcommand{\id}{{\rm{id}}}

\newcommand{\im}{{\rm{Im}}}
\newcommand{\mult}{{\rm{mult}}}

\renewcommand{\ker}{\mathop{{\rm{Ker}}}\nolimits}

\newcommand{\cl}{{\rm{closure}}}

\newcommand{\ity}{{\infty}}

\newcommand{\e}{\varepsilon}
\newcommand{\fin}{\hspace*{\fill}$\Box$\vspace*{2mm}}

\newcommand{\cG}{{\cal G}}

\newcommand{\cI}{{\cal I}}

\newcommand{\cW}{{\cal W}}

\newcommand{\cZ}{{\cal Z}}

\newcommand{\bC}{{\mathbb C}}

\newcommand{\bP}{{\mathbb P}}

\newcommand{\bZ}{{\mathbb Z}}
\newcommand{\bX}{{\mathbb X}}

\begin{document}

\title[Vanishing cycles of a meromorphic function]
{On the vanishing cycles of a meromorphic function on the complement of its 
poles}
\author{Dirk Siersma}
\address{D.S.: Mathematisch Instituut, Universiteit Utrecht, PO
Box 80010, \ 3508 TA Utrecht,
 The Netherlands.}  
\email{siersma@math.uu.nl}
\author{Mihai Tib\u ar}  
\address{M.T.:  Math\' ematiques, UMR 8524 CNRS,
Universit\'e de Lille 1, \  59655 Villeneuve d'Ascq, France.}
\email{tibar@agat.univ-lille1.fr}


\thanks{The authors benefited from the ``Research in Pairs" 
program at Oberwolfach, 
supported by Volkswagen Foundation, and from the CNRS-NWO collaboration in research program, with support from the French side for the second author.}

\begin{abstract}
We study vanishing cycles naturally attached 
to a meromorphic function with isolated singularities, in both local and global settings.
\end{abstract}
\maketitle
\setcounter{section}{0}
\section{Introduction}\label{s:intro}

 In the local setting, we consider a meromorphic function germ $p_x/q_x \colon (\cZ, x) 
\dashrightarrow \bP^1(\bC)$, where $(\cZ, x)$ is the germ of a complex manifold and $p_x$ and $q_x$ are holomorphic germs at the point $x$. 
By definition, $p'_x/q'_x$ is equal to $p_x/q_x$ if and only 
if  there exists a 
holomorphic germ $u$ such that $u(x) \not= 0$ and that $p_x=up'_x$ and $q_x= 
uq'_x$.

In the global setting, a meromorphic function $p/q \colon \cZ \dashrightarrow \bP^1(\bC)$ is defined as the ratio of two sections, $p$  
and $q$, of a holomorphic line bundle $L\to \cZ$ over a connected compact 
complex 
manifold $\cZ$. We consider only examples of projective $\cZ$, which case insures the existence of global meromorphic functions, by Kodaira embedding theorem.

 Our constructions and results in the two settings are completely similar and parallel. This is why we shall adopt in this paper a unique notation for both situations:  $p/q \colon \cZ \dashrightarrow \bP^1(\bC)$ can alternatively mean a meromorphic germ or a global meromorphic function.
 
The meromorphic function induces a holomorphic function $f: X \to \bC$, 
on the complement of $X:=\cZ\setminus \{ q=0\}$ of the pole locus $\{ q=0\}$.
Then we call {\em vanishing homology} the relative homology $H_*(X,F; \bZ)$, where $F$ denotes a general fiber of $f: X \to \bC$. 
Since $f$ is a non-proper function, vanishing cycles may appear not only because of the critical 
points of $f$, but also because of a certain  non-regular behavior of the fibers of $f$ in the neighbourhood of the pole locus $\{ q=0\}$, which is more subtle to detect.
 We interpret the later phenomenon in terms of singularities of the meromorphic function, to which we give a precise meaning. Our study concerns the class of meromorphic functions with {\em isolated singularities}, including possible singularities occuring in the indeterminacy locus.
 This is a class of non-generic pencils, far beyond the class of generic pencils that is currently considered in Lefschetz theory. We send the reader to \cite{Ti-ijm} and \cite{Ti-surv} for comments on the relations to Lefschetz theory.
 
 While the homology of $X$ and that of $F$
 might be very complicated, it turns out (Theorem \ref{t:attach})
 that the vanishing homology $H_*(X,F; \bZ)$ is concentrated in dimension $n = \dim X$ and moreover that
the space $X/F$ has the homotopy type of a bouquet of spheres 
$\bigvee S^n$.
 The proof needs new technical ingredients, due to the general notion of  
singularity that we use here.  Proposition 
\ref{p:inv} and Lemma \ref{l:contr} are crucial in this respect. 

 We prove that the {\em polar Milnor number} which we attach to an 
isolated $\cG$-singularity is the number of vanishing cycles 
which are ``concentrated" at this singularity (Proposition 
\ref{t:conc} and Corollary \ref{c:num}). Alternatively, the polar Milnor number  $\lambda_\xi$ at some point $\xi$  coincides 
with the jump of the usual local Milnor number of the pencil of hypersurfaces.
Vanishing cycles can now be detected by a multibranch jump formula (Theorem 
\ref{t:numbers}).

We also show that the global vanishing homology decomposes into a 
direct 
sum of vanishing homologies at the atypical fibres, with localization in case of $\cG$-isolated 
singularities (Proposition \ref{t:conc} and Corollary \ref{c:num}).
 We prove a general global 
Picard type formula for the monodromy (Proposition 
\ref{p:piclef}). Significant examples and particular cases are 
treated in \S \ref{s:ex}. 

The vanishing homology in the context of meromorphic functions extends the one of local proper holomorphic functions, 
which has been initiated by Milnor \cite{Mi} and developed in many aspects ever since. 
 The meromorphic setting also generalizes the study of the topology of 
polynomial functions via their 
singularities at infinity and the study of 
one-parameter families of non compact hypersurfaces, developped in the last 
decade. (We send the reader to \cite{Ti-surv} for more details on these connections.) New phenomena may occur: unlike the polynomial case, where 
singularities
at infinity ``stay at the same place'' for all fibers,
 in the meromorphic case
singularities can split or even disappear (see Example \ref{ex:1bis}).

  The results addressed here are based on our 1999 preprint \cite{ST-x}, which 
  has not been published. Ever since, 
there appeared several papers which quote it, like  \cite{DN2}, 
\cite{GLM2}, \cite{Ti-surv}. The recent survey \cite{Ti-surv} treats aspects 
of the topology of meromorphic functions on singular spaces, reviewing some of the results presented in the preprint \cite{ST-x} and in this paper. Connectivity results of Lefschetz type via nongeneric pencils (i.e. global meromorphic functions) are proved in \cite{Ti-ijm}. 
 From a quite different point of view, meromorphic germs are discussed in 
[GLM1,2], where the main interest is 
 the zeta-function of the monodromy. Classifications aspects have been explored by Arnold \cite{Ar}, who found the list of simple germs of 
meromorphic functions under natural equivalence relations.

%
\section{Singularities of $f$ along the indeterminacy locus}

The definitions in this section naturally extend
the ones used in case of polynomials and certain classes of regular functions on affine varieties, see \cite{Ti-t, Ti-surv} and also \cite{Ti-ijm}.
According to our convention in the Introduction, we treat in parallel the local and global settings, using a common notation.
 
\begin{definition}
 We call {\em completed space} the global hypersurface $\bX$ of 
$\cZ\times \bP^1$ defined
 by:
 \[  sp(z)-tq(z)=0 , \]
 where $[t:s]\in \bP^1$. Then $\bX$ is the analytic closure in $\cZ\times \bP^1$ 
of the graph of $f 
:X\to\bC$.
 We denote by $\pi : \bX \to \bP^1$ the projection on $\bP^1$. This is a proper holomorphic
function which 
extends $f \colon X\to \bC$.  The space $\bX^{\pol} := \bX \cap  
(\{ p=q=0\}\times \bP^1 )$ is a divisor of 
$\bX$ and $X$ is identified with  $\bX\setminus (\bX^{\pol}\cup \bX_\ity)$, where $\bX_\ity := \pi^{-1}(\{ s=0\})$.  
\end{definition}


We endow $\bX$ with a 
locally finite Whitney stratification $\cW$ such that $X$ is a 
stratum. In case of a germ at $x\in \cZ$ of a meromorphic function  
one considers the germ 
of such a Whitney stratification at the line $\{x\}\times \bP^1\subset \bX$.
In both situations, the local 
finiteness of the strata implies, by using Thom-Mather Isotopy Lemma \cite{T} \` a la Verdier \cite{Ve}, the following finiteness result. 
 
\begin{proposition}\label{p:1}
The stratified projection $\pi :\bX \to \bP^1$ with respect to 
$\cW$ is locally topologically trivial
over $\bP^1\setminus \Lambda_f$, for some finite set $\Lambda_f\subset\bP^1$.
In particular, the restriction $\pi_{| X} = f$ is a locally trivial {\rm 
C}$^\ity$-fibration over $\bC\setminus \Lambda_f$.
\fin 
\end{proposition}

We shall call the {\em set of atypical values}, and denote it 
by 
$\Lambda_f$, the 
minimal set $\Lambda_f$ which satisfies Proposition \ref{p:1}.
 Notice that, if we take  
the meromorphic 
function $q/p$ instead of $p/q$, then we get the same space $\bX$ (even if
$X$ will be $Z\setminus \{ p=0\}$ instead of $Z\setminus \{ q=0\}$). In 
particular, $\Lambda_{\frac{1}{f}}$ can be identified with  $\Lambda_f$
by the isomorphism $[t:s] \mapsto [s:t]$.

For any subset $A\subset \bP^1$, we denote $\bX_A := \pi^{-1}(A)$, 
$F_A := f^{-1}(A)$ and 
the general fibre $F := F_t = X\cap 
f^{-1}(t)$, for some $t\not\in \Lambda_f$.
Let $D$ be a small 
disc centered at $a\in \Lambda_f$, such that $D\cap \Lambda_f = \{ a\}$.

A crucial problem in investigating the topology of the fibres of $f$ 
is how to detect and to control the change of topology. This is an open problem, 
in general, even for 
holomorphic germs or polynomial functions, but well understood in case the singularities are 
isolated. 
  

\begin{definition}
Let $( \cZ, x) \dashrightarrow \bP^1$ be a germ of a meromorphic 
function. To every $a\in \bP^1$, one associates the germ $\pi : (\bX , (x, a)) \to \bP^1$. By 
restriction to $X = \bX \setminus (\bX^{\pol}\cup \bX_\ity)$, this defines the 
germ $f_|: (X , (x, a)) \to \bC$, which we denote by $f_{x,a}$. 
\end{definition}

For some fixed $x\in 
\{ q=p=0\}$, one has a {\em one-parameter family of 
germs} $(\bX , (x, a))$, indexed over $a\in \bP^1\setminus \{ s=0\}$.
Unlike the case of holomorphic germs, here the point  
$(x, a)$ is not in $X$ but in its closure.
So the germ $f_{x,a}$ is uniquely 
determined by the 
determination of the point $a = [p(x): q(x)]$.

Fourtheron, take a Whitney stratification $\cW$ of $\bX$ which has $X$ as 
open stratum, remarking that $\Sing\bX \subset \bX^\pol$. For all small enough radii $\e$ of a ball $B_\e (x,a)\subset \cZ\times \bP^1$ centered at $(x,a)$, the sphere $S_\e = 
\partial \bar B_\e(x,a)$ intersects transversally all the 
finitely many  strata in the neighbourhood of $(x,a)$. This defines a {\em Milnor-L\^e fibration} (cf \cite{Le-oslo}), i.e. a locally trivial fibration $\pi : \bX_{D^*} \cap B_\e(x,a)  \to D^*$ over a small enough punctured disk $D^*\subset \bP^1$ centered at $a$,  which restricts 
to a locally trivial fibration on the complement of $\bX^{\pol}$, namely:
\begin{equation}\label{eq:2}
f_| : F_{D^*} \cap B_\e(x,a)  \to D^*.
\end{equation} 

This fibration will be called the 
 Milnor-L\^e fibration of the function germ $f_{x,a}$ at the 
point $(x, a) \in \bX$. It depends on the point $(x,a)$. In particular, the radius $\e$ of the ball depends on the point $a$.\footnote{One may compare to [GLM1,2], where different definitions have been used.}
From Proposition \ref{p:1} it follows that, since $\pi$ is 
stratified-transversal to $\bX$ over $\bP^1 \setminus \Lambda_f$, the 
fibration  
$f_| : F_D \cap B_\e (x,a) \to D$ is a trivial 
fibration, for all but a finite number of values of $a\in \bP^1$, 
where $x\in \{ p=q=0\}$ is fixed.

We now endow $\bX$ with a ``partial Thom stratification", cf \cite{Ti-t}. 
  Suppose that $\bX$ is endowed with a complex stratification
$\cG = \{ \cG_\alpha\}_{\alpha \in S}$ such that $X$ is a stratum and 
$\bX^{\pol}$ is a union of strata. If  $\cG_\alpha \cap 
\overline{\cG_\beta} \not=
\emptyset$ then, by definition, $\cG_\alpha \subset 
\overline{\cG_\beta}$ and in
this case we  write $\cG_\alpha < \cG_\beta$.

Let $(x,a)\in \bX^{\pol}\setminus \bX_\ity$. We consider on $\bX$ the 
{\em Thom} (a$_q$) {\em regularity condition} at $(x,a)$, see e.g. 
\cite[ch. I]{GWPL} for the definition. 
In terms of 
the {\em relative conormal} (see \cite{Te}, \cite{HMS}), 
the condition (a$_q$) at $\xi :=(x,a)$ for the strata $\cG_\alpha$ 
and 
$\cG_\beta$
translates to the inclusion: $T_{\cG_\alpha}^* \supset (T^*_{q | 
\overline{ \cG_\beta}})_\xi$.
It is known that this condition is independent on $q$, up to 
multiplication by 
a unit, see e.g. \cite[Prop. 3.2]{Ti-t}. We therefore may and shall 
refer 
to 
this as {\em Thom regularity condition relative to $\bX^{\pol}$, at 
$(x,a)$}.

\begin{definition}\label{d:partial} 
We say that $\cG$ is a $\partial \tau$-stratification ({\em partial 
Thom stratification}) relative to $\bX^{\pol}$ if at any point $\xi\in \bX$, any 
two strata 
$\cG_{\alpha} <\cG_{\beta}$ with
$\xi \in \cG_{\alpha}
\subset \bX^{\pol}$ and $\cG_{\beta} \subset \bX\setminus \bX^{\ity}$ 
satisfy 
the Thom regularity condition relative to $\bX^{\pol}$.
\end{definition}


The Whitney stratification $\cW$ of $\bX$ considered in Proposition \ref{p:1} is an example of $\partial 
\tau$-stratification relative to $\bX^{\pol}$. This follows from 
  \cite[Th\'eor\`eme 
4.2.1]{BMM} or \cite[Theorem 3.9]{Ti-t}.
 Nevertheless the $\partial \tau$-stratifications are less demanding 
than 
Whitney stratifications and than Thom stratifications. 
\hyphenation{stra-ti-fi-cation}
  

\begin{definition}\label{d:sing}
We consider the singular locus of $\pi$ with respect to some 
$\partial \tau$-stratification relative to $\bX^{\pol}$, denoted by $\cG$ 
and we say that the following closed subset of $\bX \setminus \bX_\ity$:
\[ \Sing_{\cG}f := (\bX \setminus \bX_\ity) \cap \cup_{
 \cG_\alpha \in \cG }
 \cl(\Sing \pi_{|\cG_\alpha} )\]
 is ``the singular locus of $f$" with respect to $\cG$.
We say that $f$ has {\em isolated singularities} with respect to $\cG$ 
if 
$\dim \Sing_\cG f \le 0$. We say that $f$ has isolated singularities at 
$a\in \bC$ (or at the fibre $\bX_a$) if $\dim \bX_a\cap \Sing_\cG f \le 
0$.
\end{definition}

 The space $X$ is nonsingular and consists of one stratum. The set 
$X\cap \Sing_\cG f$ of $\cG$-singularities on $X$ is just 
the usual singular set $\Sing f\subset \cZ\setminus \{ q=0\}$. 
The singularities of the new type are $\bX^{\pol} \cap  \Sing_\cG f$.
 
 We show that singularities of this type are manageable when they are isolated. In this case 
one may localize the variation of topology of fibres, which phenomenon has been observed before
in the case of holomorphic germs, by Milnor \cite{Mi}, and in case of polynomial and regular functions 
\cite{ST}, \cite{Ti-t}. Actually the proof for 
meromorphic functions follows the 
arguments of \cite[Theorem 4.3]{Ti-t} and we leave it to the reader.

For instance, when $\dim \cZ =2$, the pencil $p/q$ has isolated $\cG$-singularities, relative to the 
coarsest partial Thom 
stratification $\cG$, if and only if the fibres of $f$ are reduced.  
\begin{proposition}\label{t:loc} 
Let $f$ have isolated singularities with respect to $\cG$ at  $a\in 
\bC$ and let $\bX_a \cap \Sing_\cG f =\{ a_1, \ldots, a_k\}$.
Then the variation of topology of the fibres of $f$ at $F_a$ is
localizable at the points $a_i$.
\fin
\end{proposition}
The localization result implies that the vanishing cycles are 
concentrated at the isolated singularities, as follows:

\begin{proposition}\label{t:conc}
Let $f$ have isolated singularities with respect to $\cG$ at $a\in \bC$ 
and let $\bX_a \cap \Sing_\cG f =\{ a_1, \ldots, a_k\}$. Let $D\subset 
\bC$ be a small enough closed disc centered at $a$ and let $s\in 
\partial 
D$. Then, for any small enough balls $B_i\subset \cZ\times \bC$ centered 
at $a_i$, we have:
\begin{enumerate}
\rm \item \it $H_*(F_D, F_s) \simeq \oplus_{i=1}^k \tilde 
H^{2n-1-*}(B_i\cap \bX_s)$.
\rm \item \it  $H_*(B_i\cap F_D , B_i\cap F_s) \simeq 
\tilde H^{2n-1-*}(B_i\cap \bX_s)$, $\forall i \in \{ 1, \ldots , k\}$.
\end{enumerate}
\end{proposition}

\begin{proof} Note first that in the local setting we have just one singular point, i.e. $k=1$.\\ 
(a). A general Lefschetz duality result (see e.g. \cite[Prop. 5.2]{Br}) 
says that, since we work with triangulable spaces, we have:
\[ H_*(F_D, F_s)\simeq H^{2n-*}(\bX_D,\bX_s).\]
 Next, the cohomology group splits, through excision, into local 
contributions, by  Proposition \ref{t:loc}:
 \[ H^{*}(\bX_D , \bX_s) = \oplus_{i=1}^k H^*(B_i\cap \bX_D, B_i\cap 
\bX_s) = \oplus_{i=1}^k \tilde H^{*-1}(B_i\cap \bX_s) , \]
 where the second equality holds because $B_i\cap \bX_s$ is 
contractible, 
for small enough ball $B_i$.
 
 (b). The same Lefschetz duality result may be applied locally to yield:
 \[ H_*(B_i\cap F_D,B_i\cap  F_s) \simeq H^{2n-*}(B_i\cap \bX_D, B_i\cap 
\bX_s).\]
Note that the decomposition $H_*(F_D, F_s)\simeq \oplus_{i=1}^k 
H_*(B_i\cap F_D, B_i\cap F_s)$ also follows from  
\ref{t:loc}. 
\end{proof}
 Our main result concerning vanishing cycles at the level of homotopy type is   the following. 

\begin{theorem}\label{t:attach} 
 Let  $f$ have isolated singularities with respect to
some $\partial \tau$-stratification $\cG$ relative to $\bX^{\pol}$. Let 
$F$ be a general fibre of $f$ and let $D\subset \bC$ be a small enough open disc centered at some $a\in \Lambda_f$. 
Then the space $X$, resp. $F_D$, is obtained from 
$F$ to
which one attaches a number of cells of real dimension $n = \dim \cZ$.
In particular we have the following homotopy equivalences:
\begin{enumerate}
\rm \item \it $X/F \simeq \bigvee S^{n}$, in the global setting.
\rm \item \it $F_D/F \simeq \bigvee S^{n}$, in the local or global setting.
\end{enumerate} 
\end{theorem}
In the local setting, by $F_D$ and $F$ we mean the intersections with some small sphere $B_\e$, so (b) should read: $B_\e \cap F_D/B_\e \cap F \simeq \bigvee S^{n}$.

We shall give the proof in \S \ref{s:polar}, after introducing a few 
technical ingredients.
The number of spheres will be discussed in Corollary \ref{c:num}.
In the local setting, Theorem \ref{t:attach} extends Milnor's bouquet theorem \cite{Mi}, whereas in the global setting, it extends the bouquet result for
polynomial functions \cite[Theorem 3.1]{ST} and is of similar flavor as \cite[Theorem 4.6]{Ti-t}.

\section{Polar curves and Milnor numbers at the indeterminacy locus}
\label{s:polar}
 We show first that an isolated $\cG$-singularity at a point of 
$\bX^\pol$ is detectable by the presence of a certain local polar locus, 
which we define in the following. 

\begin{definition}\label{d:polar}
 Let $\xi= (x,a)\in \bX^\pol \setminus \bX_\ity$ and consider a 
small neighbourhood $V\subset \cZ$ of $x$. Let $\Sing 
f$, respectively $\Sing(f, q)$, denote the singular locus of the restriction 
$f : X\cap (V\times \bC) \to \bC$, respectively $(f,q): X\cap (V\times \bC) \to \bC^2$. 
 
 The {\em polar locus} $\Gamma_\xi(f,q)$ is the germ at $\xi$ of the 
analytic space: 
 \[ \cl \{ \Sing(f, q) \setminus (\Sing f \cup \bX^\pol)\} \subset 
\bX .\]
 \end{definition}

From the definition we get the isomorphism
$\Gamma_\xi(f,q) \simeq \Gamma_\xi(p,q)$.
  The polar locus depends on the multiplicative unit $u$, i.e. 
$\Gamma_\xi(f,qu)$ is different from $\Gamma_\xi(f,q)$, 
 meanwhile we shall prove that
it induces well defined local invariants.

\begin{proposition}\label{p:inv}
Let $\xi= (x,a)\in \bX^\pol \setminus \bX_\ity$. 
Let $f$ have an isolated $\cG$-singularity at $\xi$. Then:
\begin{enumerate}
\rm \item \it  The polar locus 
$\Gamma_\xi(f,qu)$ is either void or $\dim \Gamma_\xi(f,qu) =1$, and this does not depend on the multiplicative unit $u$.
\rm \item \it  The intersection multiplicity $\mult_\xi 
(\Gamma_\xi(f,qu), \bX_a )$ is independent on the  
unit $u$.
\end{enumerate}
\end{proposition}
\begin{proof}  
For (a). we only give the rough idea and send for details to
\cite[Prop. 4.2]{ST-x} and \cite[Prop. 3.4]{Ti-surv}.
The first claim follows by usual arguments, as in \cite{Ti-t}.
The independence on $u$ is a consequence of the independence of the
relative conormal $\bP T^*_{qu}$ proved in \cite[Prop. 3.2]{Ti-t}.

 If the polar locus is void, the 
multiplicity in (b) is zero. Suppose next that $\Gamma_\xi(f, q)$ has dimension $1$. Consider a small enough ball $B\subset 
\cZ \times \bC$ centered at $\xi$, to fit in the Milnor-L\^e fibration 
 of the function $\pi$ at $\xi$:
\begin{equation}\label{eq:le}
\pi_| : B \cap \bX_{D^*} \to D^* , 
\end{equation}
where $D\subset \bC$ is centered at $a$.  
 The notation $\Gamma(f, q)$ will stay for the representative in $B$ 
of 
the germ  $\Gamma_\xi(f, q)$.
 We may choose $D$ so small, that for all $s\in \partial D$, those 
intersection points $\bX_s \cap \Gamma(f, q)$ which tend to $\xi$ when 
$s$ tends to  $a$, are inside $B$. This is possible because $\Gamma(f, q)$ 
is a curve which cuts $\bX^\pol$ at $\xi$.

We shall compute the homology $H_*(B\cap \bX_s)$ of the Milnor fibre of 
the fibration (\ref{eq:le}). Inside $B$, the restriction  of the 
function 
$q$ to $B\cap \bX_s$ has a finite number of isolated singularities, 
which 
are precisely the points of intersection $B\cap \bX_s \cap \Gamma(f, 
q)$. 

We claim that the space $B\cap \bX_s \cap  
q^{-1}(\hat\delta)$ is  
contractible, for small enough disc $\hat\delta \subset \bC$ centered at 
$0$. To prove it, we need the following:

\begin{lemma}\label{l:contr}
Let $f$ have isolated $\cG$-singularities at $\xi$.  Let $B$ be a small 
enough ball at $\xi$ such that the 
sphere $S:= \partial \bar B$ cuts transversely all those finitely many 
strata 
of $\cG$ which have $\xi$ in 
their closure and does not intersect other strata.

 Then, there exist small enough discs  $D$ and $\delta$ such that 
$(\pi,q)^{-1}(\nu)$ is transverse to $S$, for all $\nu\in D\times 
\delta^*$.
\end{lemma}
\begin{proof}

If the statement was not true, then there would exist a sequence of
points $\eta_i \in S \cap (\bX\setminus \bX^\pol)$ tending to a point
$\eta\in S \cap \bX_a \cap \bX^\pol$, such that
the intersection of tangent spaces
$T_{\eta_i}f^{-1}(f(\eta_i)) \cap T_{\eta_i}q^{-1}(q(\eta_i))$ is 
contained in $T_{\eta_i}(S\cap X)$.
Assuming, without loss of generality, that the following limits exist, 
we 
get the inclusion:
\begin{equation}\label{eq:limit}
\lim T_{\eta_i}f^{-1}(f(\eta_i)) \cap \lim 
T_{\eta_i}q^{-1}(q(\eta_i)) \subset
\lim T_{\eta_i}(S\cap X) .
\end{equation}
Let $\cG_\alpha \subset \bX^\pol$ be the
stratum containing $\eta$. Remark that $\dim \cG_\alpha \ge 2$, since 
$\bar \cG_\alpha \ni \xi$ and $\pi \pitchfork_\eta \cG_\alpha$. This 
implies that $\dim \cG_\alpha \cap \bX_a \ge 1$.

We have, by the definition of the stratification $\cG$, that $\lim 
T_{\eta_i} q^{-1}(q(\eta_i)) \supset T_\eta \cG_\alpha$ and obviously 
$T_\eta (\cG_\alpha \cap \bX_a) \subset T_\eta \cG_\alpha$.
On the other hand, $\lim T_{\eta_i}f^{-1}(f(\eta_i))\supset T_\eta
(\cG_\alpha \cap \bX_a)$, since $\pi \pitchfork_\eta 
\cG_\alpha$. In conclusion, the intersection in (\ref{eq:limit}) 
contains 
$T_\eta
(\cG_\alpha \cap \bX_a)$. But, since $S \pitchfork_\eta \cG_\alpha$, 
the limit $\lim T_{\eta_i}(S\cap X)$ cannot contain $T_\eta
(\cG_\alpha \cap \bX_a)$ and this gives a contradiction.
\end{proof}
 

Let $\hat\delta$ be so small that $B\cap \bX_s \cap 
q^{-1}(\hat\delta) \cap \Gamma(f, q) = \emptyset$. By the Lemma 
\ref{l:contr} above and by 
choosing appropriate $D$ and $\hat\delta$, the map $q : B\cap \bX_s \cap 
q^{-1}(\hat\delta^*) \to \hat\delta^*$ is a locally trivial fibration.
 Therefore $B\cap \bX_s \cap q^{-1}(\hat\delta)$ is homotopy equivalent, 
by 
retraction, to the central fibre $B\cap \bX_s \cap 
\bX^\pol$.  This proves our claim.

We now remark that the central fibre   
$B\cap \bX_s \cap \bX^\pol$ is just the complex link at $\xi$ of 
the space $\bX^\pol$. The space $\bX^\pol$ is a product $(\{ q=0\} \cap 
\{ p=0\}) \times \bC$ at $\xi$, along the projection axis $\bC$, hence its complex link is contractible, and so is $B\cap \bX_s \cap 
q^{-1}(\hat\delta)$.

Pursuing the proof of Proposition \ref{p:inv}, we observe that $B\cap 
\bX_s$ 
is homotopy equivalent to $B\cap \bX_s \cap q^{-1}(\delta)$, for $D$ and 
$\delta$ 
like 
in Lemma \ref{l:contr} and, in addition, the radius of $D$ much 
smaller 
than the radius of $\delta$. This supplementary condition is meant to 
insure 
that $\Gamma(f, q) \cap B \cap \bX_s = \Gamma(f, q) \cap B \cap 
\bX_s 
\cap q^{-1}(\delta)$.  

Now, the total space $B\cap \bX_s \cap q^{-1}(\delta)$ is built up by 
attaching 
to the space $B\cap \bX_s \cap  q^{-1}(\hat\delta)$, which is 
contractible,  
a finite number 
of cells of dimension $n-1$, which correspond to the Milnor numbers of 
the 
isolated singularities of the function $q$ on $B\cap \bX_s \cap q^{-1} 
(\delta \setminus \hat\delta)$. The sum of these numbers is, by 
definition, 
the 
intersection multiplicity 
$\mult_\xi(\Gamma(f, q) , \bX_a)$.

We have proven that:
 \begin{equation}\label{eq:mult} 
 \dim H_{n-1}(B\cap \bX_s) = \mult_\xi ( \Gamma(f , q), \bX_a) \ \ 
\mbox{and} \ \ \tilde H_i(B\cap \bX_s) = 0, \ \mbox{for} \ i\not= n-1. 
 \end{equation}
 When replacing  all over in our proof the function $q$ by $qu$, we get 
the same relation (\ref{eq:mult}), with $qu$ instead of $q$. This 
concludes our proof of \ref{p:inv}.
\end{proof} 

The above proof shows that $B\cap \bX_s$ is, homotopically, a ball to 
which one attaches a certain number of $(n-1)$-cells. We therefore get:

\begin{corollary}\label{c:spheres}
 Let $f$ have an isolated $\cG$-singularity at $\xi$. The  fibre $B\cap 
\bX_s$ of the local fibration \rm (\ref{eq:le}) \it is homotopy 
equivalent to a bouquet of spheres $\bigvee S^{n-1}$.
 \fin
\end{corollary}
 

\begin{definition}\label{d:lambda}
We denote the number of spheres by $\lambda_\xi := \dim H_{n-1} (B\cap 
\bX_s)$ and call it the {\em polar Milnor number} at $\xi$.
We say that $f$ {\em has vanishing cycles at $\xi$} if 
$\lambda_\xi 
>0$.

In the global setting, then we denote by 
$\lambda_a$ the sum of the polar Milnor numbers at singularities on 
$\bX_a 
\cap \bX^\pol$ and also denote $\lambda = \sum_{a\in \Lambda_f} \lambda_a$.
\end{definition}
 

\subsection{Proof of Theorem \ref{t:attach}}
We take back the notations of Theorem \ref{t:attach}. Since $\Sing_\cG 
f\subset \bX \setminus \bX_\ity$ 
is a finite set of points, the 
variation of topology of the fibres of $f$ is
localizable at those points (cf. Proposition \ref{t:loc}). Let $\bX_a 
\cap \Sing_\cG f = \{ a_1, \ldots , a_k\}$, with $k=1$ in the local setting.  

For some point $a_i\in X\cap \Sing_\cG f$, by Milnor's  classical result for holomorphic functions with isolated singularity 
\cite{Mi}, it follows that the pair $(B_{\e,i} \cap F_{D_a}, 
B_{\e,i} \cap F_s)$ is ($n -1$)-connected, 
where $s \in D_a^*$.

In case $a_i \in \bX^{\pol}\cap \Sing_\cG f$, we may 
invoke the following lemma, which is a version of 
a result by Hamm and L\^e \cite[Corollary 4.2.2]{HL} for our partial Thom 
stratification:
\begin{lemma}\label{l:hammle} \rm (\cite[Cor. 2.7]{Ti-t}) \it 
The pair $(B_{\e,i} \cap F_{D_a}, B_{\e,i} \cap
F_s)$ is $(n -1)$-connected, where $s \in D_a^*$.
\fin
\end{lemma}

We conclude that the space $X$ is built up starting from a fibre 
$F$, then
moving it within a fibration with a finite number of isolated
singularities. By the above connectivity result and by Switzer's result
\cite[Proposition 6.13]{Sw}, at each singular point one has
to attach a number of cells of dimensions $\ge n$. In fact the cells to be attached are of dimension precisely $n$, by the following reason.
We may apply the duality Proposition \ref{t:conc}(b) and invoke Corollary \ref{c:spheres}, which show that the relative homology 
 $H_*(F_{D_a},F_s)$ is concentrated in dimension $n$.
 
  Then one can map a bouquet of $n$ spheres into $F_D/F_s$
such that this map is an isomorphism in homology. This implies, by Whitehead's theorem, that the map induces an isomorphism of homotopy groups. (Remark that $F_D/F_s$ is simply connected whenever $n\ge 2$). Since we work with analytic objects, therefore triangulable, the space $F_D/F_s$ is a CW-complex. For CW-complexes, weak homotopy equivalence coincides with homotopy equivalence. 

Let us remark that the 
total number of $n$-cells is the sum of the local Milnor numbers, resp. 
the polar Milnor numbers.
This ends the proof of our theorem.

 As consequence, we get the relative Betti numbers 
 (see Proposition \ref{t:conc}). This may be compared to similar formulas in case of polynomial functions \cite{ST}. 

\begin{corollary}\label{c:num}
 Let $f$ have isolated $\cG$-singularities at $a\in \bC$ with respect to 
some $\partial \tau$ stratification $\cG$. Then $H_{j}(F_D, F_s) = 0$  for $j\not= n$ and: 
\[  b_{n}(F_D, F_s) = (-1)^n \chi(F_D, F_s) = \mu_a + \lambda_a  , 
\]
where $\mu_a$ is the sum of the Milnor numbers of the singularities of 
$F_a$ and $\lambda_a$ denotes the sum of the polar Milnor numbers at 
$\bX_a\cap \bX^\pol$.

In particular, if $f$ has isolated $\cG$-singularities at all fibres, 
then:
\[  b_{n}(X,F) = (-1)^n \chi(X, F) = \mu + \lambda \ , \hspace*{2cm} \ 
H_{j}(X, F) = 
0, \ \mbox{for} \ j\not= n ,\]
where $\mu$ is the total Milnor number of the singularities of 
$f$ on $\cZ\setminus \{ q=0\}$ and $\lambda$ is the total polar Milnor 
number at $\bX^\pol\setminus \bX_\ity$.
\fin
\end{corollary}

\section{Vanishing cycles in special cases and examples}\label{s:ex}
The singular locus 
$\Sing \bX \subset \cZ\times \bC$ is contained in $\bX^\pol$ and can be 
complicated. We have $\Sing \bX \setminus \bX_\ity = \cup_{t\in 
\bC}(\Sing \bX_t) \cap \bX^\pol$.
However, $\bX^\pol \setminus \Sing \bX$ is a  Whitney stratum and 
$\Sing \bX$ is a union of Whitney strata, in the canonical Whitney 
stratification $\cW$ of $\bX$ which has $X$ as a stratum.

We shall consider here a $\partial \tau$-stratification $\cG$ which may be
coarser than $\cW$ (and which exists, by Definition \ref{d:partial} and 
the remark following it). Then $\Sing_\cG f \cap \bX^\pol \subset \Sing 
\bX$. Indeed, this follows from the fact that the space $\bX^\pol 
\setminus \bX_\ity$ is locally a product $\{ q=p=0\} \times \bC$ and the 
projection $\pi$ is transversal to it off $\Sing \bX$.

In particular, for $n=2$, $f$ has isolated $\cG$-singularities at $a\in \bC$ if 
and only if $F_a$ is reduced.

 Let $\xi =(x, a)\in \bX^\pol \setminus \bX_\ity$.
We assume in the following that $\dim_\xi \Sing \bX_a = 0$. 
 This implies that $\dim_\xi \Sing_\cG f \le 0$ and that the germ 
$(\Sing \bX, \xi)$ is either a curve or just the point $\xi$.
 If it is a curve, then it can have several branches and its intersection with 
$\bX_s$ is, 
say, $\{ \xi_1(s), \ldots 
,\xi_k(s)\}$, for any $s\in D^*$, where $D\subset \bC$ is a small enough 
disc at $a$. 

The germs $(\bX_s, \xi_i(s))$ are germs of hypersurfaces with isolated 
singularity. Let $\mu_i(s)$ denote the Milnor number of   $(\bX_s, 
\xi_i(s))$. Then $\sum_{i=1}^{k} \mu_i(s) \le \mu(a)$. Equality 
may hold 
only if $k=1$, by the well known non-splitting result of L\^e D.T 
\cite{Le-a}.
In general, we have:
\begin{theorem}\label{t:numbers}
Let $\dim_\xi \Sing \bX_a =0$ and $\dim_\xi \Sing \bX = 1$. Then:
 \[ \lambda_\xi = \mu(a) - \sum_{i=1}^{k} \mu_i(s) .\]
\end{theorem}
\begin{proof}
 The hypothesis implies that the germ of $\Sing_\cG f$ at $\xi$ is  
the point $\xi$ or it is void. 
 For any $s\in D$ small enough, the germ $(\bX_s , \xi_i(s))$ is locally 
defined by the function:
 \[ h = p-tq : (\cZ\times \bC, \xi_i(s)) \to \bC.\]
 We have that, locally at $\xi$, the singular locus $\Sing h$ is equal 
to 
 $\Sing \bX$, in particular included into $\bX^\pol$.
 Consider the map $(h,t) : (\cZ\times \bC, \xi_i(s)) \to \bC^2$. Note 
that the polar locus $\Gamma_\xi(h,t)$ is a curve or it is void, since 
$\xi$ is an isolated $\cG$-singularity.
 Following \cite{Le-oslo}, there is a fundamental 
system of privileged polydisc neighbourhoods of $\xi$ in $\cZ\times 
\bC$, 
of the form $(P_\alpha \times D'_\alpha)$, where $D'_\alpha\subset \bC$ 
is a disc at $a$ and $P_\alpha$ is a polydisc at $x\in \cZ$ such 
that the map 
 \[ (h,t) : (\cZ \times \bC) \cap (P_\alpha \times D'_\alpha) \cap 
(h,t)^{-1} (D_\alpha \times D'_\alpha) \to D_\alpha \times D'_\alpha \]
 is a locally trivial fibration over $(D_\alpha^* \times D'_\alpha)  
\setminus \im(\Gamma(h,t))$. We chose $D_\alpha$ and $D'_\alpha$ such 
that 
 $\im(\Gamma(h,t)) \cap \partial(\overline{D_\alpha^* \times D'_\alpha}) 
= \im(\Gamma(h,t)) \cap (D_\alpha^* \times 
\partial(\overline{D'_\alpha}))$.
 Let $s \in \partial \overline{D'_\alpha}$. Observe that $t^{-1}(s) \cap 
(P_\alpha \times D'_\alpha)$ is contractible, since it is the Milnor fibre 
of 
the linear function $t$ on a smooth space.
 This is obtained, up to homotopy type, by attaching to 
$(h,t)^{-1}(0,s) \cap (P_\alpha \times D'_\alpha)$ a certain number $r$ 
of $n$-cells, equal to the sum of the Milnor numbers of the function 
$h_| 
: t^{-1}(s) \cap (P_\alpha \times D'_\alpha) \to D_\alpha$.
 Since we have the homotopy equivalence $(h,t)^{-1}(0,s) \cap (P_\alpha 
\times 
D'_\alpha) \simeq B\cap \bX_s$, we get, by Corollary \ref{c:spheres} and 
Definition \ref{d:lambda}, that $r= \lambda_\xi$.
 
 Now $(h,t)^{-1}(\eta,s) \cap (P_\alpha \times D'_\alpha)$ is homotopy 
equivalent to the Milnor fibre of the germ $(\bX_a, \xi)$, which has 
Milnor number $\mu(a)$. The space $t^{-1}(s) \cap (P_\alpha \times 
D'_\alpha)$ is obtained from $(h,t)^{-1}(\eta,s) \cap (P_\alpha \times 
D'_\alpha)$ by attaching exactly $r$ cells of dimension $n$ (coming 
from 
the polar intersections) and of a number of $n$-cells coming from the 
intersections with $\Sing h$. This number of cells is, by definition,
 $\sum_{i=1}^k \mu_i(s)$. We get the equality
 $ \mu(a) = r + \sum_{i=1}^k \mu_i(s)$.
 Since $r= \lambda_\xi$, our proof is done.
\end{proof}


\begin{remark}
If in the hypothesis of Theorem \ref{t:numbers} the dimension of $\Sing 
\bX$ is not $1$ but $0$, then the result still holds, with the remark 
that in this case $\mu_i(s)=0$, $\forall i$ and $\forall s\in D^*$.
Hence $\lambda_\xi = \mu(a)$.
\end{remark}

We give in the remainder three examples, one on $\bP^2(\bC)$ with no 
singularities 
in the complement of $\bX^\pol$ and two on a nonsingular quadratic surface 
in  $\bP^3(\bC)$. 

\begin{example}\label{ex:1}
 $E^{a,b}_{p,q} :$ \ \ $\displaystyle f= \frac{x(z^{a+b} + 
x^ay^b)}{y^pz^q}$, with $a+b+1 = p+q$ and $a,b,p,q \ge 1$.\\
This defines a meromorphic function on $\bP^2(\bC)$.
 For some $t\in \bC$, the space $\bX_t$ is given by:
 \begin{equation}\label{eq:1}
 x(z^{a+b} + x^ay^b) = t y^pz^q
 \end{equation}
 We have $\bX^\pol \setminus \bX_\ity = \{ [1:0:0], [0:1:0], [0:0:1]\} 
\times \bC$. According to Theorem \ref{t:numbers}, we look for jumps in 
the Milnor number within the family of germs (\ref{eq:1}):
 \begin{enumerate}
 \item at $[1:0:0]$, chart $x=1$. No jumps, since uniform Brieskorn type 
$(b, a+b)$.
 \item at $[0:1:0]$, chart $y=1$. For $t\not= 0$, Brieskorn type $(a+1, 
q)$, with $\mu(t) = a(q-1)$. If $t=0$, then we have $x^{a+1} + xz^{a+b}=0$ 
with $\mu(0) = a^2 + ab + b$ and the jump at $\xi = ([0:1:0], 0)$ is 
$\lambda_\xi = a^2 + ab + b - a(q-1) = b +ap$, by Theorem \ref{t:numbers}.
 \item at $[0:0:1]$, chart $z=1$. No jumps, since type $A_0$ for all $t$.
 \end{enumerate}
 We get the total jump $\lambda = b+ap$. A straightforward computation 
shows that $\mu =0$. 
 
 The fibres of $f$ can be described as follows. If $t=0$, we have $c+1$ 
disjoint copies of $\bC^*$, where $c= \gcd (a,b)$; hence $\chi(F_0) =0$. 
If $t\not= 0$, we compute $\chi(F) =- (b+ap)$, by a branched covering 
argument. The vanishing homology is concentrated in dimension 2. 
 Taking $X= \bC^2 \setminus \{ y=0\}$, we get the Betti number
  $b_2(X,F) = \chi(X,F) = \chi(X) - \chi(F) = 0 + (b+ap) = b+ap$.
  It follows $b_2(X,F) = \lambda + \mu$, which agrees with Corollary 
\ref{c:num}. 
\end{example}

\begin{example}\label{ex:1bis}
Let $\cZ \subset \bP^3$ be the nonsingular 
hypersurface given by $h = x^2 + z^2 + yw=0$ and
consider the meromorphic function $y/x$. It has its axis 
$\{ x=y=0\}$ tangent to $h=0$ at $[0:0:0:1]$.  
By computations we get $\mu =0$, $\lambda =1$ (jump $A_0 \to A_1$, at 
$t=0$). The general fibre $F$ is contractible, the special fibre $F_0$ is 
$\bC \sqcup \bC$  and $X$ is homotopy equivalent to $S^2$. All the 
connected components of fibres are contractible, however the global vanishing 
homology is generated by a relative 2-cycle. We may remark here that the jump 
$A_0 
\to A_1$ cannot 
occur in case of polynomial functions at infinity.
\end{example}

\begin{example}\label{ex:2}
 Consider the meromorphic function $\displaystyle f=\frac{x(z^2 + 
xy)}{z^3}$ on the nonsingular hypersurface $\cZ \subset \bP^3$ given by $h = yw 
+ x^2 - z^2=0$.  Then $\bX^\pol \setminus \bX_\ity = [0:0:0:1] \times \bC 
\cup [0:1:0:0] \times \bC$, where $[x:y:z:w]$ are the homogeneous 
coordinates in $\bP^3$.
 
 Along $[0:0:0:1] \times \bC$, in the chart $w=1$ and coordinates $x$ and 
$z$ on $\bX$, we have the family of curves (germs of $\bX_t$):
 \begin{equation}\label{eq:i}
 x( z^2 - x^3 + xz^2) = t z^3.
\end{equation}
 For all $t$, this is a $D_5$ singularity, so no jumps.
 
 Along $[0:1:0:0] \times \bC$, in the chart $y=1$ and, again,  $x$ and $z$ 
as coordinates on $\bX$, we have the family of curves (germs of $\bX_t$):
 \begin{equation}\label{eq:ii}
 x( z^2 +x) = t z^3.
\end{equation}
 This has type $A_2$ if $t\not= 0$ and $A_3$ if  $t=0$. Thus the jump at 
$\xi :=([0:1:0:0],0)$ is $\lambda_\xi =1$ and the total jump is $\lambda = 
1$.
 
 By simple computations, we get $\mu =2$, since there are two singular 
fibres, $F_{\pm 1}$, with $A_1$-singularities. There are 3 atypical 
fibres: $F_0 \simeq \bC^* \sqcup \bC^*$, $F_{\pm 1} \simeq \bC^*$ and the 
general fibre $F \simeq \bC^{**}$. Since $X \simeq S^2$, we get $b_2(X, F) 
= 2 - (-1) = 3$ global vanishing cycles, $X/F \simeq \bigvee_3 S^2$.
 
\end{example}


\section{Monodromy fibration and a global Picard phenomenon}\label{s:mon}

We continue to consider the local and global settings in the same time. Let $\hat D \subset \bC$ be a closed disc, big enough such that $\hat D 
\supset\Lambda_f\setminus 
\infty$ , where $\ity$ denotes the point $[1:0] \in \bP^1$. Let $D_i \subset\hat 
D$ be a 
small enough closed disc at $a_i \in \Lambda_f$, such 
that $D_i\cap \Lambda_f = \{ a_i\}$. Take a point $s$ on the boundary 
of 
$\hat D$ and, for each $i$, a path $\gamma_i \subset \hat D$ from $s$ 
to 
some fixed point $s_i\in\partial D_i$, with the usual conditions: the 
path $\gamma_i$ has no self intersections and does not intersect any 
other path $\gamma_j$, except at the point $s$. By Proposition 
\ref{p:1}, 
the fibration $f: X \setminus f^{-1}(\Lambda_f) \to \bC \setminus 
\Lambda_f$ is locally trivial, hence we may use excision in the pair 
$(F_{\hat D}, F_s)$ and get an isomorphism (induced by the inclusion of 
pairs):
\begin{equation}\label{eq:iso}
 \oplus_{a_i\in\Lambda_f} H_*(F_{D_i}, F_{s_i})  \to H_*(X, F_s) ,
\end{equation}
This shows that each inclusion $(F_{D_i}, F_{s_i}) \subset 
(X,F_{s_i})$ induces
 an injection in homology $H_*(F_{D_i}, F_{s_i}) \hookrightarrow H_*(X, 
F_{s_i})$. We also get by excision the following split exact sequence:
\[
 0\to H_*(F_{D_i}, F) \to H_*(X,F) \to 
\oplus_{a_j\in\Lambda_f, j\not= i} H_*(F_{D_j}, F) \to 0.
\]

We next consider the monodromy $h_i$ around an atypical value $a_i \in 
\Lambda_f$. 
This is induced by a counterclockwise loop around the small circle 
$\partial D_i$. The monodromy acts on the pair $(X,F)$ and we denote its 
action in homology by $T_i$. 

The following sequence of maps:
\begin{equation}\label{eq:wang} 
H_{q+1}(X,F) \stackrel{\partial}{\to} H_q(F) \stackrel{w}{\to}
H_{q+1}(F_{\partial D_i},F) \stackrel{i_*}{\to} H_{q+1}(X,F),
\end{equation}

where $w$ denotes the Wang map (which is an isomorphism, by the K\"unneth 
formula), 
gives, by composition, the map
$T_i - \id : H_{q+1}(X,F) \to H_{q+1}(X,F)$.

This overlaps the first two maps in the following sequence:
$ H_q(F) \stackrel{w}{\to}
H_{q+1}(F_{\partial D_i},F) \stackrel{i_*}{\to} 
H_{q+1}(X,F)\stackrel{\partial}{\to} H_q(F)$.
The last arrow in the sequence (\ref{eq:wang}) fits in the commutative 
diagram:
\[ \begin{array}{rcl}
  H_{q+1}(F_{\partial D_i},F) & \stackrel{i_*}{\longrightarrow}  & 
H_{q+1}(X,F) \\
  \searrow & \  & \nearrow  \\
 \ &  H_{q+1}(F_{D_i}, F) & \ 
\end{array} \] 
where all three arrows are induced by inclusion. 

It follows that 
the submodule of ``anti-invariant cycles" $\cI_*(T_i) :=\im (T_i - \id : 
H_* (X, F)\to H_* (X, F))$ 
is 
contained in the direct summand $H_*(F_{D_i}, F)$ of $H_* (X, F))$.
If $\cI_*$ denotes the submodule generated by $\cI_*(T_i)$, for all $a_i 
\in \Lambda_f$, then:
\begin{equation}\label{eq:sum}
\cI_* = \oplus_{a_i\in\Lambda_f} \cI_*(T_i) .
\end{equation}

 Using Picard's decomposition of the monodromy, Lefschetz has proven the famous relation for the monodromy around a simple nodal singularity on a nonsingular ambient space, wellknown as {\em Picard-Lefschetz formula}. The following result describes a {\em global Picard phenomenon}.
\begin{proposition}\label{p:piclef}
Identify $H_*(X,F)$ to 
$\oplus_{a_i\in \Lambda_f} H_*(F_{D_i}, F)$ by the isomorphism {\rm 
(\ref{eq:iso})}. Then, for $\omega \in H_*(X,F)$, we have:
\[ T_i (\omega) = \omega + \psi_i(\omega) , \]
for some $\psi_i(\omega) \in H_*(F_{D_i}, F)$.
\fin
\end{proposition}


In the global setting, when specializing to a homologically contractible total space $X$, the natural 
$\partial$-map  $H_{*} (X, F) \to  \tilde 
H_{*-1}(F)$
becomes an isomorphism and we get:
$ \tilde H_*(F) = \oplus_{a_i\in\Lambda_f} H_{* +1}(F_{D_i}, F)$.
  This occurs for instance in case of a polynomial function 
$g : \bC^n \to \bC$, for which $X = \bC^n$.  
Results on invariant cocycles were obtained in 
\cite{NN}. In our more general setting, these results can also be proved by dualizing from homology to cohomology. 
One obtains in this way statements about invariant cocycles $\ker (T^i -\id : H^*(X,F) \to H^*(X,F))$ instead of anti-invariant cycles. 

  In the particular case 
of polynomial functions, the above Picard formula (extracted from our preprint \cite{ST-x}) was independently noticed in \cite{NN2} and \cite{DN1}. 
 
 We end by an easy consequence, remarked
in the special case of polynomial functions in \cite{DN2}, which holds in our more general setting of local and global
meromorphic functions.

\begin{corollary}\label{c:pis}
 Assume that the number of paths is $l$ and the paths $\gamma_1, \ldots 
\gamma_l$ are counterclockwise ordered. Then the Coxeter element $T_\ity 
:=  T_l \circ \cdots \circ T_1$ determines the generators $T_i$, $\forall 
i\in\{ 1, \ldots , l\}$.
\end{corollary}
\begin{proof}
Use the direct sum decomposition (\ref{eq:sum}) together with
 the following adapted decomposition of  
$T_\ity - \id$, where we couple two-by-two consecutive terms:
$T_\ity - \id = (T_l \circ \cdots \circ T_1 - T_{l-1} \circ \cdots \circ 
T_1) + \ldots + (T_1 - \id)$.
\end{proof}

%

\end{document}